\documentclass[reqno]{amsart}
\usepackage{amssymb}

\title[Ambrosetti-Prodi problem with degenerate potential]
{Ambrosetti-Prodi problem with degenerate potential and Neumann boundary condition}

\author{Du\v{s}an D. Repov\v{s}}

\address{Faculty of Education and Faculty of Mathematics and Physics,
University of Ljublja\-na \& Institute of Mathematics, Physics and Mechanics, SI-1000 Ljubljana, Slovenia}\email{dusan.repovs@guest.arnes.si}

\newtheorem{theorem}{Theorem}[section]

\def\hunu{H^1(\Omega;|x|^\alpha)}

\def\RR{{\mathbb R}}

\def\phi{\varphi}
\def\div{\mathop{\rm div}\,}
\def\diva{\div (|x|^\alpha\nabla u)}

\def\phi{\varphi}

\def\di{\displaystyle}
\def\ri{\rightarrow}
\def\intom{\int_\Omega}

\def\bb{\begin{equation}}
\def\bbb{\end{equation}}

\begin{document}


\begin{abstract}
We study the degenerate elliptic equation $-\diva =f(u)+t\phi(x)+h(x)$ in a bounded open set $\Omega$ with homogeneous Neumann boundary condition, where $\alpha\in(0,2)$ and $f$ has a linear growth. The main result establishes the existence of real numbers $t_*$ and $t^*$ such that the problem
 has at least two solutions if $t\leq t_*$, there is
 at least one solution if $t_*<t\leq t^*$, and no solution exists
for all $t>t^*$. The proof combines {\it a priori} estimates with topological degree arguments.
\end{abstract}
\maketitle

\section{Introduction}
Let $\Omega\subset\RR^N$ be a bounded open set with smooth boundary. In their seminal paper \cite{ap}, A.~Ambrosetti and G.~Prodi studied the following semilinear elliptic problem
\bb\label{pr1}
\left\{\begin{array}{lll}
&\di \Delta u+f(u)=v(x)&\quad\mbox{in}\ \Omega\\
&\di u=0&\quad\mbox{on}\ \partial\Omega,
\end{array}\right.
\bbb
where the nonlinearity $f$ is a function whose derivative crosses the first (principal) eigenvalue $\lambda_1$ of the Laplace operator in $H^1_0(\Omega)$, in the sense that
 $$0<\lim_{t\ri -\infty}\frac{f(t)}{t}<\lambda_1<\lim_{t\ri +\infty}\frac{f(t)}{t}<\lambda_2.$$
 By using the abstract approach developed in \cite{ap}, A.~Ambrosetti and G.~Prodi have been able to describe the exact number of solutions of problem \eqref{pr1} in terms of $v$, provided that $f''>0$ in $\RR$. More precisely, they proved that there exists a closed connected manifold $A_1\subset C^{0,\alpha}(\overline\Omega)$ of codimension 1 such that $C^{0,\alpha}(\overline\Omega)\setminus A_1=A_0\cup A_2$ and problem \eqref{pr1} has exactly zero, one or two solutions according as $v$ is in $A_0$, $A_1$, or $A_2$. The proof of this pioneering result is based upon an extension of Cacciopoli's mapping theorem to some singular case.

A cartesian representation of $A_1$ is due to M.S. Berger and E. Podolak \cite{berger}, who observed that it is convenient to write problem \eqref{pr1} in an equivalent way, as follows. Let
$$Lu:=\Delta u+\lambda_1u,\qquad g(u):=f(u)-\lambda_1u$$
and
$$v(x):=t\phi(x)+h(x)\quad\mbox{with}\quad \intom h(x)\phi(x)dx=0.$$
In such a way, problem \eqref{pr2} is equivalent to
\bb\label{pr2}
\left\{\begin{array}{lll}
&\di Lu+g(u)=t\phi(x)+h(x)&\quad\mbox{in}\ \Omega\\
&\di u=0&\quad\mbox{on}\ \partial\Omega,
\end{array}\right.
\bbb
with $g''>0$ in $\RR$ and
$$-\lambda_1<\lim_{t\ri -\infty}\frac{g(t)}{t}<0<\lim_{t\ri +\infty}\frac{g(t)}{t}<\lambda_2-\lambda_1.$$

Under these assumptions, M.S. Berger and E. Podolak \cite{berger} proved that there exists $t_1$ such that problem \eqref{pr2} has exactly zero, one or two solutions according as $t<t_1$, $t=t_1$, or $t>t_1$. The proof of this result is based on a global Lyapunov-Schmidt reduction method.

For related developments on Ambrosetti-Prodi problems we refer to H.~Amann and P.~Hess \cite{amann}, D.~Arcoya and D.~Ruiz \cite{arcoya}, E.B.~Dancer \cite{dancer}, P.~Hess \cite{hess}, J.L.~Kazdan and F.W.~Warner \cite{kazdan}, J.~Mawhin \cite{maw1, maw2}, etc.

The present paper is concerned with the Ambrosetti-Prodi problem in relationship with the contributions of P.~Caldiroli and R.~Musina \cite{caldi}, who initiated the study of Dirichlet elliptic problems driven by the differential operator ${\rm div}\, (|x|^\alpha\nabla u)$, where $\alpha\in (0,2)$. This operator is a model for equations of the type
\bb\label{bao}-{\rm div}\, (a(x)\nabla u)=f(x,u)\quad x\in\Omega,\bbb
where the weight $a$ is a non-negative measurable function that is allowed to have ``essential" zeros at some points or even to be unbounded. According to
R.~Dautray and J.-L.~Lions \cite[p. 79]{dautray}, equations like \eqref{bao} are introduced as models for several
physical phenomena related to equilibrium of anisotropic continuous media
which possibly are somewhere ``perfect" insulators or ``perfect" conductors. We also refer to the works by M.K.V.~Murthy and G.~Stampacchia
[16] and by M.S.~Baouendi and C.~Goulaouic \cite{bg}, concerning degenerate elliptic operators (regularity of solutions and spectral
theory). Problem \eqref{bao} also has some interest in the framework of optimization
and $G$-convergence, cf. B.~Franchi, R.~Serapioni, and F.~Serra Cassano \cite{franchi}.
For degenerate phenomena in nonlinear PDEs we also refer to 
G.~Fragnelli and D.~Mugnai \cite{frag} and M. Nursultanov and G. Rozenblum \cite{nurs}.

The present paper deals with the study of the Ambrosetti-Prodi problem in the framework of the degenerate elliptic operator studied in \cite{caldi}. A feature of this work is that the analysis is developed in the framework of Neumann boundary conditions.

\section{Main result and abstract setting}
Let $\alpha\in (0,2)$ and let $\Omega\subset\RR^N$ be a bounded open set with smooth boundary.

Consider the following nonlinear problem
\bb\label{1}
\left\{\begin{array}{lll}
&\di -\diva =f(u)+t\phi(x)+h(x)&\quad\mbox{in}\ \Omega\\
&\di \frac{\partial u}{\partial\nu}=0&\quad\mbox{on}\ \partial\Omega.
\end{array}\right.
\bbb
We assume that $f:\RR\ri\RR$ is a continuous function such that
\bb\label{f}\limsup_{t\ri -\infty}\frac{f(t)}{t}<0<\liminf_{t\ri +\infty}\frac{f(t)}{t}\bbb
and there exists $C_f>0$ such that
\bb\label{f1}|f(t)|\leq C_f(1+|t|) \ \mbox{for all $t\in\RR$.} \bbb

Since the first eigenvalue of the Laplace operator with respect to the Neumann boundary condition is zero, condition \eqref{f} asserts that the nonlinear term $f$ crosses this eigenvalue.

Next, we assume that $\phi$, $h\in L^\infty(\Omega)$ and
\bb\label{phi}\phi\geq 0,\ \phi\not\equiv 0\ \mbox{in}\ \Omega.\bbb

Since $\alpha>0$, the weight $|x|^\alpha$ breaks the invariance under translations and can give rise to an abundance of existence results, according to the geometry of the open set $\Omega$.

For $\zeta\in C^\infty_c(\Omega)$ we define
$$\|\zeta\|^2_\alpha:=\intom (|x|^\alpha|\nabla\zeta|^2+\zeta^2)dx$$
 and we consider the function space
 $$\hunu:=\mbox{closure of $C^\infty_c(\overline{\Omega})$ with respect to the $\|\,\cdot\,\|_\alpha$-norm}.$$

 It follows that $\hunu$ is a Hilbert space with respect to the scalar product
 $$\langle u,v\rangle_\alpha:=\intom (|x|^\alpha\nabla u\cdot\nabla v+uv)dx,\quad\mbox{for all $u.v\in\hunu$}.$$

 Moreover, by the Caffarelli-Kohn-Nirenberg inequality (see \cite[Lemma 1.2]{caldi}), the space $\hunu$ is continuously embedded in $L^{2_\alpha^*}(\Omega)$, where $2_\alpha^*$ denotes the corresponding critical Sobolev exponent, that is, $2_\alpha^*=2N/(N-2+\alpha)$.

 We say that $u$ is a solution of problem \eqref{1} if $u\in\hunu$ and for all $v\in\hunu$
 $$\intom |x|^\alpha\nabla u\cdot\nabla vdx=\intom f(u)vdx+t\intom \phi vdx+\intom hvdx.$$

   Since the operator $Lu:=-\diva$ is uniformly elliptic on any strict subdomain $\omega$ of $\Omega$ with $0\notin\overline\omega$, the standard regularity theory can be applied in $\omega$. Hence, a solution $u\in\hunu$ of problem \eqref{1} is of class $C^\infty$ on $\Omega\setminus\{0\}$.
   We refer to H.~Brezis \cite[Theorem IX.26]{brezis} for more details.

The main result of this paper extends to the degenerate setting formulated in problem \eqref{1} the abstract approach developed by P.~Hess \cite{hess} and F.~de Paiva and M.~Montenegro \cite{depaiva}. For related properties on Ambrosetti-Prodi problems with Neumann boundary condition, we refer to  A.E.~Presoto and F.~de Paiva \cite{presoto}, E.~Sovrano \cite{sovrano}, A.~V\'elez-Santiago \cite{velez, velez1}.

\begin{theorem}\label{th1}
Assume that hypotheses \eqref{f}, \eqref{f1} and \eqref{phi} are fulfilled. Then there exist real numbers $t_*$ and $t^*$ with $t_*\leq t^*$ such that the following properties hold:

(a) problem \eqref{1} has at least two solutions solution, provided that $t\leq t_*$;

(b) problem \eqref{1} has at least one solution, provided that $t_*<t\leq t^*$;

(c) problem \eqref{1} has no solution, provided that $t>t^*$.
\end{theorem}

\subsection{Strategy of the proof}
Let $C_f$ be the positive constant defined in hypothesis \eqref{f1} above, and assume that $v\in L^2(\Omega)$. Consider the following linear Neumann problem:
\bb\label{2}
\left\{\begin{array}{lll}
&\di -\div(|x|^\alpha\nabla w)+C_fw =v&\quad\mbox{in}\ \Omega\\
&\di \frac{\partial w}{\partial\nu}=0&\quad\mbox{on}\ \partial\Omega.
\end{array}\right.
\bbb
With the same arguments as in \cite[Chapter IX, Exemple 4]{brezis}, problem \eqref{2} has a unique solution $w\in\hunu$. This defines a linear map
$$L^2(\Omega)\ni v\mapsto w\in\hunu.$$

It follows that the linear operator $T:L^\infty(\Omega)\ri \hunu$ defined by $Tv:=w$ is compact. We also point out that if $v\geq 0$ then $w\geq 0$, hence $T$ is a positive operator.

We observe that $u$ is a solution of problem \eqref{1} if and only if $u$ is a fixed point of the nonlinear operator $$S_t(v):=T(f(v)+C_fv+t\phi +h).$$

Thus, solving problem \eqref{1} reduces to finding the critical points of $S_t$.

\section{Proof of the main result}
We split the proof into  steps (\ref{s1})-(\ref{s4}).

\subsection{Problem \eqref{1} has no solutions for large $t$}\label{s1} In fact, we show that a necessary condition for the existence of solutions of problem \eqref{1} is that the parameter $t$ should be small enough.

We first observe that hypothesis \eqref{f} implies that there are positive constants $C_1$ and $C_2$ such that
$$f(t)\geq C_1|t|-C_2\quad\mbox{for all $t\in\RR$}.$$

Assuming that $u$ is a solution of problem \eqref{1}, we obtain by integration
$$\begin{array}{ll}
0&\di =\intom f(u)dx+t\intom\phi dx+\intom hdx\\
&\di\geq C_1\intom |u|dx-C_2|\Omega|+t\intom\phi dx+\intom hdx\\
&\di \geq -C_2|\Omega|+t\intom\phi dx+\intom hdx.\end{array}$$
It follows that a necessary condition for the existence of solutions of problem \eqref{1} is
$$t\leq\frac{\di C_2|\Omega|-\intom hdx}{\di\intom\phi dx}\,.$$

\subsection{Problem \eqref{1} has solutions for small $t$: The preliminary step}\label{s2}
In this subsection, we prove that for any $\rho>0$ there exists $t_\rho\in\RR$ such that for all $t\leq t_\rho$ and all $s\in [0,1]$ we have
\bb\label{3} v\not=sS_t(v)\quad\mbox{for all $v\in L^\infty(\Omega)$, $\|v^+\|_\infty=\rho$}.\bbb

Our argument is by contradiction. Thus, there exist three sequences $(s_n)\subset [0,1]$, $(t_n)\subset\RR$ and $(v_n)\subset L^\infty(\Omega)$ such that
$\lim_{n\ri\infty}t_n=-\infty$, $\|v_n^+\|_\infty=\rho$ and
\bb\label{4}v_n=s_nS_{t_n}(v_n)\quad\mbox{for all}\ n\geq 1.\bbb

By hypothesis \eqref{f1} we have
\bb\label{5}\begin{array}{ll}
f(v_n)+C_fv_n&\di\leq C_f+C_f|v_n|+C_fv_n\\
&\di= C_f+2C_fv_n^+\leq C_f+2C_f\rho.\end{array}\bbb

Using the definition of $S$ and the fact that $T$ is a positive operator, relations \eqref{4} and \eqref{5} yield
$$\begin{array}{ll}
v_n&\di =s_nS_{t_n}(v_n)=s_nT(f(v_n)+C_fv_n+t_n\phi +h)\\
&\di\leq s_nT(C_f+2C_f\rho+t_n\phi +h),\end{array}$$
hence
$$v_n^+\leq s_n[T(C_f+2C_f\rho+t_n\phi +h)]^+.$$

Let
$$w_n:=C_f+2C_f\rho+t_n\phi +h.$$
It follows that $w_n$ is the unique solution of the problem
$$\left\{\begin{array}{lll}
&\di -\div (|x|^\alpha w_n)+C_fw_n =C_f+2C_f\rho+t_n\phi +h&\quad\mbox{in}\ \Omega\\
&\di \frac{\partial w_n}{\partial\nu}=0&\quad\mbox{on}\ \partial\Omega.
\end{array}\right.$$
Dividing by $t_n$ (recall that $\lim_{n\ri\infty}t_n=-\infty$) we obtain
$$\left\{\begin{array}{lll}
&\di -\div \left(|x|^\alpha \frac{w_n}{t_n}\right)+C_f\,\frac{w_n}{t_n} =\phi +\frac{C_f+2C_f\rho +h}{t_n}&\quad\mbox{in}\ \Omega\\
&\di \frac{\partial}{\partial\nu}\left(\frac{w_n}{t_n} \right)=0&\quad\mbox{on}\ \partial\Omega.
\end{array}\right.$$

But
$$\lim_{n\ri\infty}\frac{C_f+2C_f\rho +h}{t_n}=0.$$
So, by elliptic regularity (see \cite[Theorem IX.26]{brezis}),
$$\frac{w_n}{t_n}\ri T\phi\quad\mbox{in}\ C^{1,\beta}(\overline\Omega\setminus\{0\})\ \mbox{as}\ n\ri\infty.$$
Next, by the strong maximum principle, we have
$$T\phi >0\quad\mbox{in}\ \Omega$$
and
$$\frac{\partial T\phi}{\partial\nu}(x)<0\quad\mbox{for all} \ x\in\partial\Omega\ \mbox{with}\ T\phi (x)=0.$$

We deduce that for all $n$ sufficiently large
$$\frac{w_n}{t_n}>0\quad\mbox{in}\ \Omega,$$
which forces
$$w_n^+=0\quad\mbox{for all $n$ large enough.}$$
But $$v_n^+\leq s_nw_n^+\leq w_n^+,$$
hence
$$\rho=\|v_n^+\|_\infty\leq\|w_n^+\|_\infty =0,$$
 a contradiction.
 This shows that our claim \eqref{3} is true.

\subsection{Problem \eqref{1} has solutions for small $t$: The intermediary step}\label{s3}
In this subsection, we prove that for any $t\in\RR$ there exists $\rho_t>0$ such that for all $s\in [0,1]$ we have
\bb\label{33} v\not=sS_t(v)\quad\mbox{for all $v\in L^\infty(\Omega)$, $\|v^-\|_\infty=\rho_t$}.\bbb

Fix arbitrarily $t\in\RR$. Assume that there exist $s\in [0,1]$ and a function $v$ (depending on $s$) such that
$v=sS_t(v)$. It follows that $v$ is the unique solution of the problem
\bb\label{34}
\left\{\begin{array}{lll}
&\di -\div(|x|^\alpha\nabla v)+C_fv =s(f(v)+C_fv+t\phi +h)&\quad\mbox{in}\ \Omega\\
&\di \frac{\partial v}{\partial\nu}=0&\quad\mbox{on}\ \partial\Omega.
\end{array}\right.
\bbb

By hypotheses \eqref{f} and \eqref{f1}, there exist positive constants $C_3$ and $C_4$ with $C_3<C_f$ such that
\bb\label{ctrei}f(t)\geq -C_3t-C_4\quad\mbox{for all $t\in\RR$}.\bbb
Returning to \eqref{34} we deduce that
$$\begin{array}{ll}
-\div(|x|^\alpha\nabla v)+C_fv&\di\geq s(-C_3v-C_4+C_fv+t\phi+h)\\
&\di= s[(C_f-C_3)v+t\phi+h-C_4].\end{array}$$
Therefore
\bb\label{35}
\left\{\begin{array}{lll}
&\di -\div(|x|^\alpha\nabla v)+[sC_3+(1-s)C_f]v \geq s(t\phi +h-C_4)&\quad\mbox{in}\ \Omega\\
&\di \frac{\partial v}{\partial\nu}=0&\quad\mbox{on}\ \partial\Omega,
\end{array}\right.
\bbb
where $$0<C_3\leq sC_3+(1-s)C_f\leq C_f.$$

Let $w$ denote the unique solution of the Neumann problem
\bb\label{36}
\left\{\begin{array}{lll}
&\di -\div(|x|^\alpha\nabla w)+[sC_3+(1-s)C_f]w = s(t\phi +h-C_4)&\quad\mbox{in}\ \Omega\\
&\di \frac{\partial w}{\partial\nu}=0&\quad\mbox{on}\ \partial\Omega.
\end{array}\right.\bbb
By \eqref{35}, \eqref{36} and the maximum principle, we deduce that
\bb\label{37}w\leq v\quad\mbox{in}\ \Omega.\bbb

Moreover, since $C_3\leq C_3\leq sC_3+(1-s)C_f\leq C_f$ for all $s\in[0,1]$, we deduce that the solutions $w=w(s)$ of problem \eqref{36} are uniformly bounded. Thus, there exists $C_0=C_0(t)>0$ such that
\bb\label{38}\|w\|_\infty\leq C_0\quad\mbox{for all}\ s\in[0,1].\bbb

Next, relation \eqref{37} yields
$$v^-=\max\{-v,0\}\leq\max\{-w,0\}=w^-\quad\mbox{in}\ \Omega.$$
Using now the uniform bound established in \eqref{38}, we conclude that our claim \eqref{33} follows if we choose $\rho_t=C_0+1$.

\subsection{Problem \eqref{1} has solutions for small $t$: The final step}\label{s4}
Let $\rho>0$ and let $t_\rho$ be as defined in subsection \ref{s2} such that relation \eqref{3} holds. We prove that problem \eqref{1} has at least one solution, provided that $t\leq t_\rho$.

Fix $t\leq t_\rho$ and let $\rho_t$ be the positive number defined in subsection \ref{s3}. Consider  the open set
$$G=G_t:=\{v\in L^\infty(\Omega);\ \|v^+\|_\infty<\rho,\ \|v^-\|_\infty<\rho_t\}.$$

It follows that
$$v\not= sS_t(v)\quad\mbox{for all $v\in\partial G$, all $s\in[0,1]$}.$$
So, we can apply the homotopy invariance property of the topological degree, see Z.~Denkowski, S.~Mig\'orski and N.S.~Papageorgiou \cite[Theorem 2.2.12]{denko}. It follows that
$${\rm deg}\, (I-S_t,G,0)={\rm deg}\, (I,G,0)=1.$$
We conclude that $S_t$ has at least one fixed point for all $t\leq t_\rho$, hence problem \eqref{1} has at least one solution.

\subsection{Proof of Theorem \ref{th1}}\label{s5}
We first show that problem \eqref{1} has a subsolution for all $t$. Fix a positive real number $t$. By \eqref{ctrei}, we have
$$f(u)+t\phi +h\geq -C_3u-C_4-|t|\,\|\phi\|_\infty-\|h\|_\infty\quad\mbox{for all}\ u\in\RR.$$
It follows that the function
$$u\equiv -\frac{|t|\,\|\phi\|_\infty+\|h\|_\infty+C_4}{C_3}$$
is a subsolution of problem \eqref{1}.

Next, with the same arguments as in the proof of Lemma 2.1 in \cite{depaiva}, we obtain that if $t$ belongs to a bounded interval $I$ then the set of corresponding solutions of problem \eqref{1} is uniformly bounded in $L^\infty(\Omega)$. Thus, there exists $C=C(I)>0$ such that for every solution of \eqref{1} corresponding to some $t\in I$ we have $\|u\|_\infty\leq C$. Since weak solutions of problem \eqref{1} are bounded, the nonlinear regularity theory of G.~Lieberman \cite{lie} implies that for every $\omega\subset\subset\Omega$ with $0\notin\overline\omega$, the set of all solutions corresponding to $I$ is bounded in $C^{1,\beta}(\overline\omega)$.

We already know (subsection \ref{s1}) that problem \eqref{1} does not have any solution for large values of $t$ and solutions exist if $t$ is small enough (section \ref{s4}). Let
$${\mathcal S}:=\{t\in\RR;\ \mbox{problem \eqref{1} has a solution}\}.$$
It follows that ${\mathcal S}\not=\emptyset$.

Let
$$t^*:=\sup{\mathcal S}<+\infty.$$

We prove in what follows that problem \eqref{1} has a solution if $t=t^*$. Indeed, by the definition of $t^*$, there is an increasing sequence $(t_n)\subset{\mathcal S}$ that converges to $t^*$. Let $u_n$ be a solution of \eqref{1} corresponding to $t=t_n$. Since $(t_n)$ is a bounded sequence, we deduce that the sequence $(u_n)$ is bounded in  $C^{1,\beta}(\overline\omega)$ for all $\omega\subset\subset\Omega$ with $0\notin\overline\omega$. By the Arzela-Ascoli theorem, the sequence $(u_n)$ is convergent to some $u_*$ in $C^{1}(\overline\omega)$, which is a solution of problem \eqref{1} for $t=t^*$.

Fix arbitrarily $t_0<t^*$. We  prove that problem \eqref{1} has a solution for $t=t_0$. We already know that problem \eqref{1} considered for $t=t_0$ has a subsolution $\underline{U}_{t_0}$. Let $u_{t^*}$ denote the solution of problem \eqref{1} for $t=t^*$. Then $u_{t^*}$ is a supersolution of problem \eqref{1} for $t=t_0$. Since $\underline{U}_{t_0}$ (which is a constant) can be chosen even smaller, it follows that we can assume that
$$\underline{U}_{t_0}\leq u_{t^*}\quad\mbox{in}\ \Omega.$$
Using the method of lower and upper solutions, we conclude that problem \eqref{1} has at least one solution for $t=t_0$.

Returning to subsection \ref{s4}, we know that for all $\rho>0$ there exists a real number $t_\rho$ such that problem \eqref{1} has at least one solution, provided that $t\leq t_\rho$. Let
$$t_*:=\sup\{t_\rho;\ \rho>0\}.$$

We already know that \eqref{1} has at least one solution for all $t<t_*$. We show that, in fact, problem \eqref{1} has at least two solutions, provided that $t<t_*$.

Fix $t_0<t_*$ and let $\rho_{t_0}$ be the positive number defined in subsection \ref{s3}. Consider  the bounded open set
$$G_{t_0}:=\{v\in L^\infty(\Omega);\ \|v^+\|_\infty<\rho,\ \|v^-\|_\infty<\rho_{t_0}\}.$$

Since $G_{t_0}$ is bounded, we can assume that $$\overline{G}_{t_0}\subset \{u\in L^\infty(\Omega); \ \|u\|_\infty<R\}=:B(0,R),$$
for some $R>0$.

Recall that if $t$ belongs to a bounded interval $I$ then the set of corresponding solutions of problem \eqref{1} is uniformly bounded in $L^\infty(\Omega)$. So, choosing eventually a bigger $R$, we can assume that
$\|u\|_\infty<R$ for any solution of problem \eqref{1} corresponding to $t\in [t_0,t^*+1]$.

Since problem \eqref{1} does not have any solution for $t=t^*+1$, it follows that
$${\rm deg}\,(I-S_{t^*+1},B(0,R),0)=0.$$ So, using
 the homotopy invariance property of the topological degree we obtain
$${\rm deg}\,(I-S_{t_0},B(0,R),0)={\rm deg}\,(I-S_{t^*+1},B(0,R),0)=0.$$
Next, using the excision property of the topological degree (see Z.~Denkowski, S.~Mig\'orski and N.S.~Papageorgiou \cite[Proposition 2.2.19]{denko}) we have
$${\rm deg}\,(I-S_{t_0},B(0,R)\setminus G_{t_0},0)={\rm deg}\,(I-S_{t^*+1},B(0,R)\setminus G_{t_0},0)=-1.$$

We conclude that problem \eqref{1} has at least two solutions for all $t<t^*$. \qed

\subsection*{Epilogue} The result established in the present paper can be extended if problem \eqref{1} is driven by degenerate operators of the type ${\rm div}\, (a(x)\nabla u)$, where $a$ is a measurable and non-negative weight in $\Omega$, which can have at most a finite number of (essential) zeros. Such a behavior holds if there exists an exponent $\alpha\in (0,2)$ such that $a$ decreases more slowly than $|x-z|^\alpha$ near every point $z\in a^{-1}\{0\}$. According to  P.~Caldiroli and R.~Musina \cite{caldi1}, such an hypothesis can be formulated as follows: $a\in L^1(\Omega)$ and there exists $\alpha\in (0,2)$ such that
$$\liminf_{x\ri z}|x-z|^{-\alpha}a(x)>0\quad\mbox{for every}\ z\in\overline\Omega.$$

Under this assumption, the weight $a$ could be nonsmooth, as the Taylor
expansion formula can easily show. For example, the function $a$ cannot be of class
$C^2$  and it cannot have bounded derivatives if $\alpha\in (0, 1)$. As established in \cite[Lemma 2.2, Remark 2.3]{caldi1}, a function $a$ satisfying the above hypothesis has a finite number of zeros in $\overline\Omega$. Notice that in such we can allow degeneracy also at some point of
its boundary.

To the best of our knowledge, no results are known for degenerate ``double-phase" Ambrosetti-Prodi problems, namely for equations driven by differential operators like
\bb\label{marc1}{\rm div}\, (|x|^\alpha \nabla u)+{\rm div}\, (|x|^\beta |\nabla u|^{p-2}\nabla u)\bbb
or
\bb\label{marc2}{\rm div}\, (|x|^\alpha \nabla u)+{\rm div}\, (|x|^\beta\log(e+|x|) |\nabla u|^{p-2}\nabla u),\bbb
where $\alpha\not=\beta$ are positive numbers and $1<p\not=2$. 

Problems of this type correspond to ``double-phase variational integrals" studied by G.~Mingione {\it et al.} \cite{mingi1,mingi2}. The cases covered by the differential operators defined in \eqref{marc1} and \eqref{marc2} correspond to a degenerate behavior both at the origin and on the zero set of the gradient. That is why it is natural to study what happens if the associated integrands are modified in such a way that, also if $|\nabla u|$ is small, there exists an imbalance between the two terms of the corresponding integrand.

\medskip
 \indent {\bf Acknowledgements.}  This research was supported by the Slovenian Research Agency program P1-0292 and grants N1-0064, J1-8131, and J1-7025.

\end{document}